# A FREE BOUNDARY MODEL FOR OXYGEN DIFFUSION IN A SPHERICAL MEDIUM

A. M. GONZÁLEZ

*Departamento de Matemática, Facultad de Ciencias Exactas, Físico-Químicas y Naturales,*
*Universidad Nacional de Río Cuarto, Ruta 36 Km 601,*
*X5804BYA Río Cuarto, Córdoba, Argentina.*
agonzalez@exa.unrc.edu.ar

J. C. REGINATO

*Departamento de Física, Facultad de Ciencias Exactas, Físico-Químicas y Naturales,*
*Universidad Nacional de Río Cuarto, Ruta 36 Km 601,*
*X5804BYA Río Cuarto, Córdoba, Argentina.*
jreginato@exa.unrc.edu.ar

D. A. TARZIA

*Departamento de Matemática and CONICET, FCE, Universidad Austral, Paraguay 1950,*
*S2000FZF Rosario, Santa Fe, Argentina.*
DTarzia@austral.edu.ar

The goal of this article is to find a correct approximated solution using a polynomial of sixth degree for the free boundary problem corresponding to the diffusion of oxygen in a spherical medium with simultaneous absorption at a constant rate, and to show some mistakes in previously published solutions.

*Keywords*: Oxygen Diffusion; Free Boundary Problem; Approximated Solution.

## 1. Introduction

The free boundary problems are generally associated with the processes of melting and freezing which have a latent heat-type condition at the interface connecting the velocity of the free boundary and the heat flux (See Refs. 1 and 2 for more details). However, there is another important class of free boundary problems that arises from the diffusion of a gas into an absorbing medium in which no such explicit condition is available at the free boundary; for example, oxygen diffusion in saturated soils spherical aggregates with constant, linear and nonlinear absorption (See Ref. 3 for more details). These free boundary problems become nonlinear due to the presence of a free boundary and for this reason their analytical solutions are difficult to obtain (See Ref. 4 for more details). A classical problem on diffusion of oxygen in a one-dimensional medium was studied in Ref. 5 in which oxygen was allowed to diffuse in a body tissue in the shape of a sheet that absorbs it at a constant rate. First the oxygen is allowed to diffuse into a medium, some of the oxygen is absorbed by the medium, thereby being removed from the diffusion process, and the concentration of oxygen at the surface of the medium is maintained constant. This phase of the problem continues until a steady state is reached in which the oxygen does not penetrate into the medium. The supply of oxygen is then cut off and the surface is sealed so that no oxygen passes in or out, the medium continues to absorb the available oxygen already in it and, as consequence, the boundary marking the furthest depth of penetration in the steady state starts to recede towards the sealed surface.

A problem on diffusion of oxygen in a spherical medium with simultaneous absorption at a constant rate is formulated in Refs. 6 and 7.

In Ref. 7, the author traced the moving boundary using the constrained integral method assuming third and fourth order polynomials. In each case he obtained a system of two ordinary differential equations for the moving boundary position and the concentration at the fixed boundary, their solution leads to the unknowns at each time step.

In Ref. 6 the author proposes a numerical method to trace the moving boundary and the concentration at the fixed surface. In this method, a double linear system of equations is considered. The first system is formed through applying first, second, third and fourth moments and using assumed profile for the concentration containing four unknowns functions of time. Numerical solution through a proposed scheme leads to the unknown functions. The second system is formed through applying the boundary conditions given in addition to another assumed condition, that is the concentration at x = 0 is unknown function of time. The results of the first system becomes an entry data for the second one leading to the concentration at the fixed surface x = 0.

The goals of these articles (Refs. 6 and 7) were to find an approximate solution applying a modification of one of the most important semi-analytical methods, the constrained integral method, proposed in Ref. 8 in order to solve implicit free boundary problems.

Taking into account the transformation from spherical to linear coordinates (see Ref. 9), the absorption process, in a dimensionless form, becomes the following free boundary problem (See Eqs. (1)–(5), pp. 671 in Ref. 6, and Eqs. (1)–(4), pp. 363–364 in Ref. 7 for more details):

$$u_t(x,t) = u_{xx}(x,t) - 1 \quad , \quad 0 < x < s(t) \ , \ t > 0 \qquad (1.1)$$

$$u(s(t),t) = u_x(s(t),t) = 0 \quad , \quad t > 0 \qquad (1.2)$$

$$u_x(0,t) = 0 \quad , \quad t > 0 \qquad (1.3)$$

$$u(x,0) = \frac{(1-x)^2}{2} \quad , \quad 0 \leq x \leq 1 \qquad (1.4)$$

$$s(0) = 1. \qquad (1.5)$$

The goal of the present article is twofold: first to find a correct approximated solution using a polynomial of six degree of the problem (1.1) to (1.5) (see Section 2) and second to show some mistakes in the solutions given in both papers Refs. 6 and 7 (see Sections 2 and 3).

## 2. About Ahmed's paper

Instead of the solution considered in Ref. 6 we propose an approximated solution of the system (1.1)–(1.5) assuming a polynomial of six degree in the spatial variable $x$ for the concentration $u$ in the region $0 < x < s(t)$, $t > 0$, that is

$$u(x,t) = a(t) + b(t)\left(\frac{x}{s(t)}\right)^2 + c(t)\left(\frac{x}{s(t)}\right)^4 + d(t)\left(\frac{x}{s(t)}\right)^6 \qquad (2.1)$$

where $a = a(t)$, $b = b(t)$, $c = c(t)$ and $d = d(t)$ are unknown parameters, to be determined, depending of the time $t$. Taking into account the conditions (1.2), (1.3) and an extra condition at the free boundary given by $u_{xx}(s(t),t) = 1$, $t > 0$, deduced from (1.1) and (1.2), the concentration profile (2.1) becomes

$$\begin{aligned}u(x,t) = a(t) &+ \left(-3a(t) + \frac{s^2(t)}{8}\right)\left(\frac{x}{s(t)}\right)^2 + \\ &+ \left(3a(t) - \frac{s^2(t)}{4}\right)\left(\frac{x}{s(t)}\right)^4 + \left(-a(t) + \frac{s^2(t)}{8}\right)\left(\frac{x}{s(t)}\right)^6\end{aligned} \qquad (2.2)$$

where $a = a(t)$ and $s = s(t)$ must be determined.

In order to obtain $a = a(t)$ (concentration at the fixed face $x = 0$) and $s = s(t)$ (the free boundary) as a function of time $t$ we take the zeroth and first moments of the differential equation (1.1) which are given respectively by

$$\int_0^{s(t)} u_t(x,t)\,dx = \int_0^{s(t)} (u_{xx}(x,t) - 1)\,dx \qquad (2.3)$$

and

$$\int_0^{s(t)} x\,u_t(x,t)\,dx = \int_0^{s(t)} x\,(u_{xx}(x,t) - 1)\,dx. \qquad (2.4)$$

Applying the Leibniz rule to the left hand side and integrating the right hand side we can write (2.3) as

$$\frac{d}{dt}\left(\int_0^{s(t)} u(x,t)\,dx\right) - u(s(t),t)\frac{ds(t)}{dt} = \left.(u_x(x,t) - x)\right|_0^{s(t)} \qquad (2.5)$$

and taking into account the boundary conditions (1.2) and (1.3), the concentration profile (2.2), and some elementary computations, we obtain the first differential equation

$$16\left(\frac{da(t)}{dt}s(t) + a(t)\frac{ds(t)}{dt}\right) + s^2(t)\frac{ds(t)}{dt} = -35\,s(t). \qquad (2.6)$$

Again from (2.4) we get

$$\frac{d}{dt}\left(\int_0^{s(t)} x\, u(x,t)\, dx\right) - s(t)\, u(s(t),t)\, \frac{d\, s(t)}{dt} = \left( x\, u_x(x,t) - u(x,t) - \frac{x^2}{2} \right)\Bigg|_0^{s(t)} \quad (2.7)$$

Therefore, by substitution in (2.7) of boundary conditions (1.2) and (1.3), the concentration profile (2.2), and some elementary computations, we obtain the second differential equation

$$6\left(\frac{d\, a(t)}{dt} s^2(t) + 2\, a(t)\, s(t)\, \frac{d\, s(t)}{dt}\right) + s^3(t)\, \frac{d\, s(t)}{dt} = 48\left(a(t) - \frac{s^2(t)}{2}\right). \quad (2.8)$$

Finally, solving (2.6) and (2.8) we get the following system of two ordinary differential equations

$$\frac{d\, s(t)}{dt} = \frac{3}{s(t)} \frac{128\, a(t) - 29\, s^2(t)}{48\, a(t) + 5\, s^2(t)} \quad (2.9)$$

$$\frac{d\, a(t)}{dt} = \frac{-84\, a(t)\, s^2(t) - 11\, s^4(t) - 768\, a^2(t)}{2\, s^2(t)\, (48\, a(t) + 5\, s^2(t))}. \quad (2.10)$$

As $s(0) = 1$ by the initial condition (1.5), we can obtain the numerical solution of the system (2.9)–(2.10) by standard methods if $a(0)$ is given. We know from the physics of the problem that $s'(t)$ is no positive, and therefore, from (2.9), we obtain

$$128\, a(t) - 29\, s^2(t) \le 0, \text{ for } 0 \le s(t) \le 1. \quad (2.11)$$

According to condition (2.11), we have assumed $a(0) = 0.2265625$. Graphs of the free boundary $s = s(t)$ vs. $t$, the unknown parameter $a = a(t)$ vs. $t$, and the corresponding concentration profile $u = u(x,t)$ vs. $x$, at different times $t$, are shown in the Figures 1, 2 and 3.

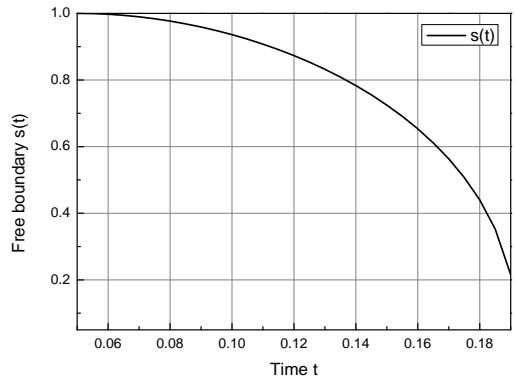

Figure 1.  The free boundary $s = s(t)$ as a function of time $t$.

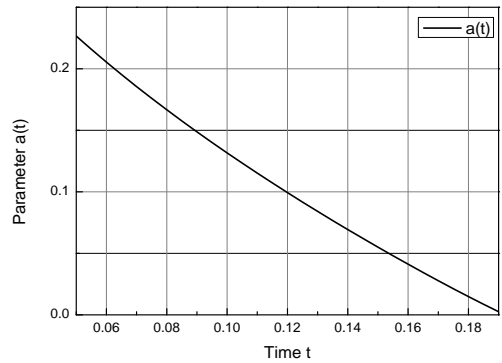

Figure 2. The parameter $a = a(t)$ as a function of time $t$.

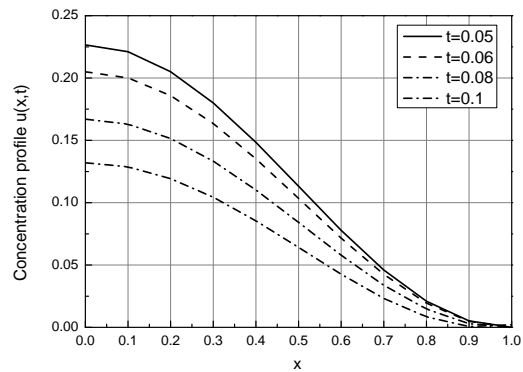

Figure 3. The concentration profile $u = u(x,t)$ versus x as a function of time $t$.

We show now in Tables 1 and 2 a comparison between our results and the one given in Ref. 10 for the concentration on the sealed surface x = 0 and the free boundary s = s(t) respectively.

Table 1. Comparative figures for the concentration at the sealed surface x = 0.

| *Time t* | *Ref. 10* | *Present* |
|---|---|---|
| 0.051 | 0.2451 | 0.227 |
| 0.060 | 0.2237 | 0.205 |
| 0.100 | 0.1423 | 0.132 |

Table 2. Comparison of position of the free boundary s = s(t).

| Time t | Ref. 10 | Present |
|---|---|---|
| 0.051 | 1.0000 | 1.000 |
| 0.060 | 0.9969 | 0.997 |
| 0.080 | 0.9756 | 0.977 |
| 0.100 | 0.9350 | 0.936 |
| 0.120 | 0.8743 | 0.873 |
| 0.140 | 0.7896 | 0.783 |
| 0.150 | 0.7356 | 0.725 |
| 0.160 | 0.6710 | 0.653 |
| 0.180 | 0.4879 | 0.440 |

From our study, given before, we remark that the equations (9) and (11) in Ref. 6 are wrong. Moreover, in that paper there exist some other differences, for example:
- from condition $u(s(t),t) = 0$ we easily obtain that $a(t) + b(t) + c(t) + d(t) = 0$ and then L.H.S. = 0 in Eq. (9) (see pp. 672 in Ref. 6),
- similar differences can be found in the use of the other moments.

With respect to the initial condition (5), given by $s(0) = 0$, in Ref. 6 we can have two cases:

Case 1: If condition (5) is correct (as it appears in the published paper) then the initial profile, given by (4), is redundant in Ref. 6. Moreover, in the numerical solution presented by the authors it is $s(0) = 1$ (See Figure 2, pp. 10 in Ref. 6). Then, the results are not consistent.

Case 2: If condition (5) is $s(0) = 1$, instead of $s(0) = 0$, then the used approximated concentration profile (6) as a polynomial in power of $((1-x)/(1-s(t))$ is not convenient. In order to solve this inconvenient we have taken, in the present paper, a polynomial in power of $(x/s(t))$.

## 3. About Çatal's paper

A better approximated solution of the system (1.1)–(1.5), with respect to the solution proposed in Ref. 7, can be found assuming a polynomial of third degree in the spatial variable *x* for the concentration *u* in the region $0 < x < s(t)$, $t > 0$, that is

$$u(x,t) = a(t) + b(t)\left(\frac{x}{s(t)}\right)^2 + c(t)\left(\frac{x}{s(t)}\right)^3 \tag{3.1}$$

where $a = a(t)$, $b = b(t)$ and $c = c(t)$ are unknown parameters, to be determined, depending of the time *t*. Taking into account the conditions (1.2) and (1.3), the concentration profile (3.1) becomes

$$u(x,t) = a(t) - 3a(t)\left(\frac{x}{s(t)}\right)^2 + 2a(t)\left(\frac{x}{s(t)}\right)^3 \qquad (3.2)$$

where $a = a(t)$ and $s = s(t)$ must be determined.

In order to obtain $a = a(t)$ and $s = s(t)$ as a function of time *t* we take the zeroth and first moments of the differential equation (1.1) which are given respectively by (2.3) and (2.4).

Applying the Leibniz rule to the left hand side and integrating the right hand side we can write (2.3) as (2.5). Then by substitution in (2.5) of the boundary conditions (1.2) and (1.3), the concentration profile (3.2), and some elementary computations, we obtain

$$\frac{d\,a(t)}{d\,t}\,s(t) + a(t)\,\frac{d\,s(t)}{d\,t} = -2\,s(t). \qquad (3.3)$$

Again from (2.4) we get (2.7). Then by substitution in (2.7) of the boundary conditions (1.2) and (1.3), the concentration profile (3.2), and some elementary computations, we obtain

$$3\left(\frac{d\,a(t)}{d\,t}\,s^2(t) + 2\,a(t)\,s(t)\,\frac{d\,s(t)}{d\,t}\right) = 20\left(a(t) - \frac{s^2(t)}{2}\right). \qquad (3.4)$$

Finally, solving (3.3) and (3.4) we get the following system of two ordinary differential equations

$$\frac{d\,s(t)}{d\,t} = \frac{4}{3}\,\frac{5\,a(t) - s^2(t)}{a(t)\,s(t)} \qquad (3.5)$$

$$\frac{d\,a(t)}{d\,t} = -\frac{2}{3} - \frac{20\,a(t)}{3\,s^2(t)}. \qquad (3.6)$$

As $s(0) = 1$ by the initial condition (1.5), we can obtain the numerical solution of the system (3.5)–(3.6) by standard methods if $a(0)$ is given. We know from the physics of the problem that $s'(t)$ is no positive, and therefore, from (3.5), we obtain

$$5\,a(t) - s^2(t) \leq 0, \text{ for } 0 \leq s(t) \leq 1. \qquad (3.7)$$

According to condition (3.7), we have assumed $a(0) = 0.2265625$. Graphs of the free boundary $s = s(t)$ vs. *t*, the unknown parameter $a = a(t)$ vs. *t*, and the corresponding concentration profile $u = u(x,t)$ vs. *x*, at different times t, are shown in the Figures 4, 5 and 6.

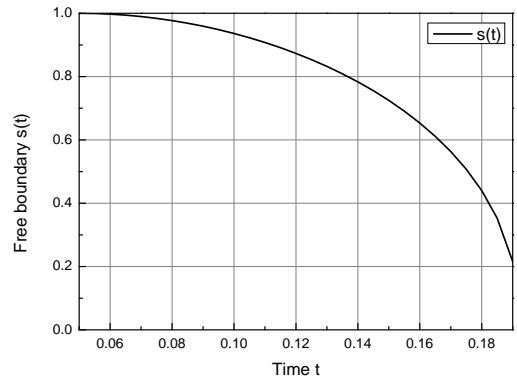

Figure 4. The free boundary $s = s(t)$ as a function of time $t$.

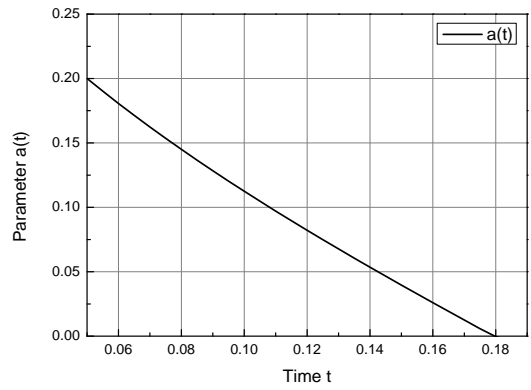

Figure 5. The parameter $a = a(t)$ as a function of time $t$.

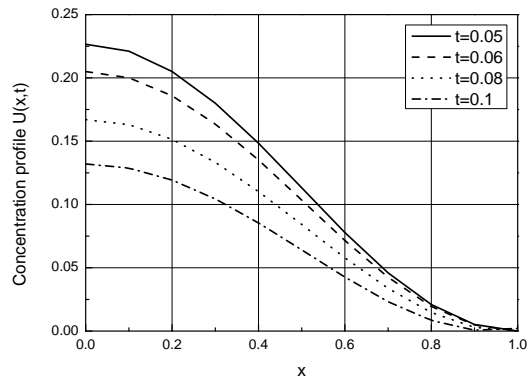

Figure 6. The concentration profile $u = u(x,t)$ versus x as a function of time $t$.

We remark that the proposed approximated polynomial of third degree in Ref. 7 is incorrect because it doesn't verify the boundary condition (1.3) and then their conclusions are not coherent.

In dimensionless form, the associated steady-state system of the evolutionary free boundary problem (1.1)-(1.5) becomes the following stationary free boundary problem (See Ref. 7, pp. 363 for more details):

$$D\, C''(X) - m = 0 \quad , \quad 0 < X < X_0 \tag{3.8}$$
$$C(X_0) = C'(X_0) = 0 \tag{3.9}$$
$$C(0) = C_0. \tag{3.10}$$

In Ref. 7, the last boundary condition (3.9) has not been proposed.

Moreover, for a correct resolution of the problem assuming a polynomial of fourth degree in the spatial variable $x$ for the concentration $u$ in the region $0 < x < s(t)$ , $t > 0$ we can see Ref. 10.

**Conclusions**

We have found a correct approximated solution to the free boundary problem (1.1) – (1.5) corresponding to the diffusion of oxygen in a spherical medium with simultaneous absorption at a constant rate by using a polynomial of sixth degree, and we show some mistakes in the solutions given by Ahmed[6] and Çatal[7].

**Acknowledgments**